\documentclass[11pt]{article}

\usepackage{amssymb,a4wide}
\usepackage{url}

\newcommand{\ndv}{v}  
\newcommand{\ndw}{w}  

\newcommand{\schm}{\mathsf{C}}

\newcommand{\bridge}[1]{\overbracket[0.6pt]{\,#1\,}} 
\newcommand{\bridget}{\bridge{\ndt}}
\newcommand{\bridgeu}{\bridge{\ndu}}
\newcommand{\bridgev}{\bridge{\ndv}}
\newcommand{\bridgew}{\bridge{\ndw}}
\newcommand{\bridgeX}{\bridge{X}}
\newcommand{\bridgeY}{\bridge{Y}}

\newcommand{\bridgeT}{\bridge{T}}
\newcommand{\bridgetree}{\bridge{\tree}}

\newcommand{\sidebyside}[2]{
    \hbox to \textwidth{
                                            \vtop{\hsize=.45 \textwidth \parindent=0pt \centering #1 \vskip1sp}
                                            \hfill
                                            \vtop{\hsize=.45 \textwidth \parindent=0pt \centering #2}
                                            } }

\newcommand{\defstyle}[1]{\textup{\textbf{#1}}}

\newcommand{\tree}{\ensuremath{\mathfrak{T}}}
\newcommand{\trt}{\mathfrak{T}}
\newcommand{\trs}{\mathfrak{S}}

\newcommand{\forestf}{\mathfrak{F}}

\newcommand{\conext}{\forestf_\tree} 
\newcommand{\conextset}{F_\tree} 
\newcommand{\fftt}{\sigma}

\newcommand{\treet}{(T; <)}

\newcommand{\nds}{s}
\newcommand{\ndt}{t}
\newcommand{\ndu}{u}

\newcommand{\tlx}[2]{{#1}^{< #2}}
\newcommand{\tleq}[2]{{#1}^{\leqslant #2}}
\newcommand{\tg}[2]{{#1}^{> #2}}
\newcommand{\tgeq}[2]{{#1}^{\geqslant #2}}

\newcommand{\pathp}{\mathsf{P}}
\newcommand{\patha}{\mathsf{A}}
\newcommand{\pathb}{\mathsf{B}}
\newcommand{\pathc}{\mathsf{C}}

\newcommand{\segm}{\mathsf{I}}
\newcommand{\segmb}{\mathsf{B}}
\newcommand{\segmj}{\mathsf{J}}

\newcommand{\stem}{\mathsf{S}}

\newcommand{\setpaths}{\mathcal{P}}

\renewcommand{\epsilon}{\varepsilon}

\newcommand{\dom}{\domain}

\newcommand{\domain}[1]{{\rm dom}(#1)}

\renewcommand{\setpaths}{\mathsf{Paths}}

\newcommand{\undl}{\mathsf{und}}
\renewcommand{\fftt}{\pi}

\newtheorem{theorem}{Theorem}

\newtheorem{corollary}[theorem]{Corollary}

\newtheorem{definition}[theorem]{Definition}
\newtheorem{example}[theorem]{Example}

\newtheorem{lemma}[theorem]{Lemma}

\newtheorem{proposition}[theorem]{Proposition}

\newenvironment{proof}[1][Proof]{\textbf{#1.} }{\ \rule{0.5em}{0.5em}}

\usepackage{mathtools}
\usepackage[font={small}]{caption}
\usepackage{xcolor}

\title
{Structural theory of trees \\
I. Branching and condensations of trees}

\author{Valentin Goranko$^1$, Ruaan Kellerman$^2$, and Alberto Zanardo$^3$  \\ 
$^1$Stockholm University, \ $^2$University of Pretoria, \ $^3$University of Padova \\ 
\textrm{valentin.goranko@philosophy.su.se}, \ 
\textrm{ruaan.kellerman@up.ac.za}, \  
\textrm{alberto.zanardo@unipd.it}}

\begin{document}

\maketitle

\begin{abstract}
Trees are partial orders in which every element has a linearly ordered set of predecessors. 
Here we initiate the exploration of the structural theory of trees with the study of different notions of \emph{branching in trees} and of \emph{condensed trees}, which are trees in which every node is a branching node.  We then introduce and investigate two different constructions of \emph{tree condensations} -- one shrinking, and the other expanding, the tree to a condensed tree.
\end{abstract}

\section{Introduction}
\label{sec:intro}
Trees are connected partial orders in which every element has a linearly ordered set of predecessors (smaller elements). After linear orders, trees are probably the most commonly and naturally occurring partial orders in a broad variety of contexts. 
Trees and tree-like structures arise not only in discrete mathematics but also in computer science (non-deterministic, concurrent, and interleaving processes, transition systems and computation trees), philosophy (models of non-determinism and branching time), game and decision theory (game trees and decision trees), and theoretical linguistics (syntax and parse trees). Trees have various applications and are associated with deep and important results in set theory (\cite{Jech}, \cite{Todorcevic}), logic (Rabin's theorem of decidability of the monadic second-order theory of some infinite trees, cf. \cite{Rabin69}), mathematics (e.g., theory of ultrametric spaces \cite{Lemin2003}, \cite{Hughes2004},  \cite{Hughes2012}), theoretical computer science (automata on trees, cf. \cite{tata2007}), graph theory and computational complexity theory (tree decompositions, structures with bounded tree width, and Courcelle's theorem, cf. \cite{DBLP:series/txcs/DowneyF13}, \cite{CourcelleEngelfriet12}).

A notable and well-explored class of trees is the class of \emph{well-founded trees}, studied mostly in a set-theoretic context, as generalizing and extending the theory of ordinals, cf. \cite{Jech}, \cite{Todorcevic}. Indeed, nearly all studies of trees found in the mathematical literature focus exclusively on that class. In fact, the well-foundedness assumption has, for historical and specific application-related reasons, been incorporated in the most commonly adopted definition of a tree in the set-theoretic tradition of their studies, and is often assumed by default.  

That is why, we have to emphasize here that in our study, trees are \emph{any connected partial orderings with linearly ordered sets of predecessors of nodes}. They are not assumed to be well-founded, and we do not think that there are any intrinsic reasons for this assumption. In fact, for most applications indicated above that assumption is not justified. Thus, we study a substantially larger class of partial orderings than well-founded trees, and we do that from a general order-theoretic perspective, not from the specific set-theoretic perspective mentioned above. We also emphasize that the well-foundedness assumption makes a very substantial difference, both in the general theory and in the particular properties of trees. Without that assumption, the study of trees remains mostly order-theoretic and extends -- in a quite non-trivial and challenging way -- the theory of linear orderings, while remaining much more specific than the general theory of partial orderings. Besides the study of well-founded trees, few general aspects of trees have been explored so far and, according to our knowledge, there has been no systematic study of the general structural theory of trees yet; certainly none coming close to the comprehensive exploration of the theory of linear orderings in \cite{Rosenstein}.

Here we initiate such a systematic exploration of the structural theory of trees. 
One general direction in the structural theory of a class $\mathcal{C}$ of mathematical structures is to identify such structures with important and desirable specific properties and then to develop and study natural generic constructions that transform or decompose any structure from $\mathcal{C}$ into one with these desirable properties. 
The present paper and the study it initiates are in that spirit. 

Three of the most important structural characteristics of trees and classes of trees are: 
\begin{enumerate}
\item The spectrum of order types of the paths (maximal chains) in the tree. 

\item The spectrum of degrees of branching at the different branching nodes in the tree.  

\item The branching structure of the entire tree. 
\end{enumerate}

Previous studies into the first of these characteristics, with emphasis on the transfer of properties and results from linear orders to trees, include \cite{GorankoZanardo}, \cite{GorankoKellerman2011}, \cite{GorankoKellerman2021}. Studies related to the degrees of branching 
include \cite{Schmerl},  \cite{Goranko1999},  \cite{SabbadinZanardo}. A general reference covering some of these results, and more, is \cite{KellermanThesis}.

Here we focus on the third of these characteristics, and begin our exploration of the general structural theory of trees with a study of the different notions of \emph{branching in trees} and of \emph{condensed trees}, which are trees in which every non-leaf node is a branching node. 
We introduce and investigate two different constructions of \emph{tree condensations}: one, \emph{shrinking} any tree into its uniquely defined `condensation', and the other -- \emph{extending} the tree to a condensed tree, containing copies of all paths in the original tree and only such paths. 
We obtain several results relating the structure of any tree to the structure of its condensations.
The notion of tree condensation extends in a certain sense the notion of condensation of linear orders defined in~\cite{Rosenstein}, which is the quotient structure obtained from a linear order by partitioning it into intervals. 
Condensations also have a topological flavour in terms of a natural notion of homeomorphism.  A further motivation to study condensed trees is that their structural theory is, generally, simpler in that both the spectra of order types of their paths and their branching structure are easier to describe than for non-condensed trees. 
Besides, both branching and condensations of trees have various applications to the general theory of trees, which go beyond the scope of the present paper, but we only mention here that both are used when axiomatising the first-order and other logical theories of some important classes of trees, and when proving the completeness of such axiomatisations, cf. \cite{KellermanThesis},  \cite{GorankoKellerman2021}.

\medskip
\textbf{Structure of the paper.}
After providing the necessary terminology and notation in Section \ref{sec:prelim}, 
we define, compare and study two notions of branching in Section \ref{sec:branching}.  We then define condensed trees and study two notions of tree condensations --  
 one shrinking, and the other expanding, the tree to a condensed tree --  
respectively in Sections \ref{sec:condensations} and \ref{sec:Condensed extensions}. We end with concluding remarks and chart our further studies in this project in Section \ref{sec:concluding}.

\section{Preliminaries}
\label{sec:prelim}

We define here some basic notions on trees, to fix notation and terminology. The reader may also consult \cite{Jech}, \cite{KellermanThesis}, and \cite{Kellerman2018} for further details.

An ordered set $\left(A;<\right)$, {with a strict partial ordering $<$},  
 is \defstyle{downward-linear} if for every $x \in A$, the set $\{ y \in A : y < x \}$ is linear;  
 it is \defstyle{downward-connected} if, for every $x,y \in A$, there exists $z \in A$ such that $z \leqslant x$ and $z \leqslant y$ (where $x \leqslant y$ is defined, as usual, as $x < y$ or $x=y$).  
 A \defstyle{forest} is a downward-linear partial order.  A \defstyle{tree} is a downward-connected forest\footnote{Note that we do not assume well-foundedness of trees, nor even the existence of a root.}. 
A  \defstyle{subtree}   
of a forest 
$\forestf = \left(F;<\right)$ 
is any substructure $\tree = \left(T;<^T\right)$ of $\forestf$ which is a tree, i.e.,  where $T$ is a non-empty downward-connected subset of $F$ and $<^T$ is the restriction of $<$ to  $T$.

The elements of a tree $\left(T;<\right)$  are called \defstyle{nodes} or \defstyle{points}. 
If a tree has a $<$-minimal node, then it is unique (by downward-connectedness)  and is called the \defstyle{root} of the tree. The $<$-maximal nodes in a tree (if there are any) are called \defstyle{leaves} of that tree. 

Next, we define various notions and notation in terms of an arbitrarily fixed tree $\tree = \treet$. 
First, for any nodes $\ndt, \ndu \in T$ we define $\ndt \smile \ndu$ to mean that $\ndt < \ndu$ or $\ndt = \ndu$ or $\ndu < \ndt$. If this holds, we say that $\ndt$ and $\ndu$ are \defstyle{comparable} nodes.
If $\ndt < \ndu$, 
the intervals $(\ndt,\ndu)$, $(\ndt,\ndu]$, $[\ndt,\ndu)$ and $[\ndt,\ndu]$ are defined as usual. 
For instance, if $\ndt < \ndu$ then $(\ndt,\ndu] := \{ x \in T : \ndt < x \leqslant u \}$, etc. 
 We also define the sets $\tlx{T}{\ndt} := \{ x \in T : x < \ndt\}$, $\tleq{T}{\ndt} := \{ x \in T : x \leqslant \ndt \}$, $\tg{T}{\ndt} := \{ x \in T : \ndt < x \}$ and $\tgeq{T}{\ndt} := \{ x \in T : \ndt \leqslant x \}$. 
 We will use analogous notation for the respective substructures (as partial orders) of the tree $\tree$ over these sets, for instance,  $\tlx{\tree}{\ndt}$ denotes $\left(\tlx{T}{\ndt};<\upharpoonright_{\tlx{T}{t}}\right)$, etc.

For  non-empty subsets $A, B \subseteq T$ we define  $A < B$ (resp. $A \leqslant  B$, $A > B$, $A \geqslant B$) when $x<y$ (resp. $x \leqslant y$, $x > y$, $x \geqslant  y$) for all $x \in A$ and $y \in B$. 
 Instead of $\{x\} < B$ we will also write $x < B$, and similarly for other relations and singleton sets.
Then, we define the sets $\tlx{T}{A} := \{ x \in T : x < A \}$ and likewise $\tleq{T}{A}$, $\tg{T}{A}$, $\tgeq{T}{A}$.  The substructures of $\trt$ that have these sets as their underlying sets will be denoted as $\trt^{<A}$, $\trt^{\leqslant A}$, $\trt^{>A}$ and $\trt^{\geqslant A}$ respectively.

More generally, given any non-empty subset $A$ of $T$, $\trt^A$ will denote the structure $\left(A;<\upharpoonright_{A}\right)$.

Note that, for any $A \not=\emptyset$, $\tlx{\tree}{A}$ and  $\tleq{\tree}{A}$ are linear orders and that $\tg{T}{A}$ and $\tgeq{T}{A}$ are empty when $A$ is not linearly ordered.  In general, if $A$ is linearly ordered then $\tg{\tree}{A}$ and $\tgeq{\tree}{A}$ are forests, while for every node $\ndt$, $\tgeq{\tree}{\ndt}$ is a tree that is rooted at $\ndt$.

A maximal linearly ordered set of nodes in a tree is called a \defstyle{path}. 
A set of nodes $\patha$ is  \defstyle{downward-closed} if $z \in \patha$ whenever $y \in \patha$ and $z < y$; respectively, $\patha$ is  \defstyle{upward-closed} if $z \in \patha$ whenever $y \in \patha$ and $y < z$.  
A non-empty linearly ordered set of nodes that is downward-closed and bounded above is called a \defstyle{stem}.  Note that a path cannot be viewed as a stem.  
A non-empty subset $\pathb$ of a path $\patha$ is called a \defstyle{branch} when it is bounded below and
upward-closed within $\patha$ (i.e.~if $x \in \pathb$ and $y \in \patha$ with $x < y$ then $y \in \pathb$).
The set of paths containing the node $\ndt$ (resp. the stem $\stem$) will  be denoted by $\setpaths_{\ndt}$ (resp. $\setpaths_{\stem})$.  

A set of nodes $\patha$ is called \defstyle{convex} if $z \in \patha$ whenever $x,y \in \patha$ and $x < z < y$.  
A convex linearly ordered set of nodes is called a \defstyle{segment}. 
A \defstyle{bridge} is a non-empty segment  $\patha$ such that, for every path $\pathp$, either $\patha \subseteq \pathp$ or $\patha \cap \pathp$ is empty. A segment $\patha$ is called a \defstyle{furcation} when it is not a bridge. 
Note that every singleton set of nodes $\{ t \}$ is a bridge.

A set $X$ of nodes is an \defstyle{antichain} if $x \not\smile y$ for all $x \not=y$ in $X$. Note that the intersection of an antichain $X$ and a linearly ordered set $Y$ of nodes is either a singleton or the empty set. The second alternative is excluded when $X$ is a maximal (by inclusion) antichain and $Y$ is a path. 
Also, note that every two distinct paths in a tree intersect in a stem.

The linear orders $\left(\mathbb{N};<\right)$, $\left(\mathbb{Z};<\right)$, $\left(\mathbb{Q};<\right)$ and $\left(\mathbb{R};<\right)$, where in each instance $<$ denotes the usual ordering of that set, will be denoted as $\omega$, $\zeta$, $\eta$ and $\lambda$ respectively.

\section{Branching in trees}
\label{sec:branching}

\subsection{Connected components and branching stems} 

\begin{definition} 
Two paths $\patha_{1}$ and $\patha_{2} $ in a tree $\tree = \treet$ are \defstyle{undivided at the stem} $\stem \subseteq \patha_{1} \cap \patha_{2}$, denoted $\patha_{1} \curlyvee_{\stem} \patha_{2}$,  
iff  $\patha_{1} \cap \patha_{2} \cap \tg{T}{\stem} \neq \emptyset$.  
Otherwise, 
$\patha_{1}$ and $\patha_{2}$ are \defstyle{branching at} $\stem$. 
If $\stem = \tleq{T}{\ndt}$, we will say that $\patha_{1}$ and $\patha_{2}$ are undivided, or branching, at $\ndt$, and $\curlyvee_{\stem}$ will  also be written $\curlyvee_{t}$.
\end{definition} 

In the sequel, some notions relative to stems are tacitly extended to nodes, like in the previous definition.

\begin{proposition} 
\label{ch2:ch2:prop:undiv}
Let $\stem$ be a stem. 
The relation $\curlyvee_{\stem}$ is an equivalence relation on $\setpaths_{\stem}$.
\end{proposition}

\begin{proof} 
Straightforward.
\end{proof}

For every path $\patha \in \setpaths_{\stem}$ we denote by $[\patha]_{\curlyvee_{\stem}}$ the equivalence class of $\patha$ with respect to 
$\curlyvee_{\stem}$.

\begin{definition}\label{def:schmerl-comp}
A \defstyle{$<$-connected component}\footnote{This definition comes from \cite{Schmerl}.}
(briefly, \defstyle{$<$-component}) of the forest $\forestf = \left(F;<\right)$  is a non-empty 
subset $\schm$ of $F$ such that: 
\begin{enumerate}
	\item
		if $\ndt \in \schm$, $\ndt' \leqslant \ndt$ and $\ndt' \leqslant \ndu$, then $\ndu \in \schm$; and
	\item
	$\schm$ is minimal (by inclusion) for the condition 1.
\end{enumerate}
\end{definition}

\begin{proposition}\label{prop:schmerl-comp}
Let  $\forestf = \left(F;<\right)$ be a forest. For any $\ndt \in F$, let
\begin{equation}\label{eq:schmerl-of-t}
\schm_\ndt = \{ \ndu : \ndt \geqslant \ndt' \leqslant \ndu , \textrm{ for some } \ndt' \in F \}\ .
\end{equation}
Then:
\begin{enumerate}
\item Every  set $\schm_\ndt$ is a $<$-component of $\forestf$, and vice-versa, every $<$-compo\-nent is of the type $\schm_\ndt$. 
\item The set of $<$-components of $\forestf$ is a partition of $F$;
\item Every $<$-component of $\forestf$ is a maximal subtree of $\forestf$.
\end{enumerate}
\end{proposition}

\begin{proof}
1. To prove that the set $\schm_\ndt$ is closed under the condition 1 in Def.~\ref{def:schmerl-comp},  
suppose $\ndt'' \in \schm_\ndt$,  $\ndt' \leqslant \ndt''$, and $\ndt' \leqslant u$.   
Then, by the definition of $\schm_\ndt$, there is a $\ndv$ such that $\ndt \geqslant  \ndv \leqslant \ndt''$. 
By downward linearity, applied to $\ndv$ and $\ndt'$, there are two cases:
 
i) $\ndv \leqslant \ndt'$.  Then, $\ndv \leqslant \ndu$, hence $\ndu \in \schm_\ndt$.

ii)  $\ndt' < \ndv$.  Then  $\ndt' < \ndt$ and  $\ndt'  \leqslant \ndu$. By definition of $\schm_\ndt$ these imply $\ndu \in \schm_\ndt$. 

To prove condition 2 (the minimality of $\schm_\ndt$), take any non-empty set $\schm  \subseteq \schm_\ndt$ which satisfies the closure condition 1 of Def.~\ref{def:schmerl-comp}. Take any $\ndu \in \schm$. Then 
$\ndu \in \schm_\ndt$, so $\ndt' \leqslant \ndu$ for some $\ndt' \leqslant \ndt$, hence 
$\ndt \in \schm$ by condition 1. Therefore, $\schm_\ndt \subseteq \schm$, again by condition 1 applied to $\schm$, hence  $\schm_\ndt = \schm$, as required. 
By the minimality, this also shows that every $<$-component is  $\schm_\ndt$ for any $\ndt$ in it.

2. It follows from the above that $\schm_\ndt \cap \schm_\ndv \not= \emptyset$ implies $\schm_\ndt = \schm_\ndv$, and hence the set of $<$-components is a partition of $F$.

3.  Clearly, every $\schm_\ndt$ is a subtree. 
For any $\ndu \not\in \schm_\ndt$, there is no $\ndw$ such that $\ndw \leqslant \ndt$  and $\ndw \leqslant \ndu$, otherwise $\ndu$ would belong to $\schm_\ndt$. Then, no subtree of $\forestf$ contains $\schm_\ndt$ properly. 
\end{proof}

\begin{corollary}\label{ch?:components of forests}
Let $\pathp$ be a path in a forest $\forestf$ and let $T_\pathp$  be the union of all paths $\patha$ such that $\pathp \cap \patha \not= \emptyset$. Then $T_\pathp$ is a $<$-component of $\forestf$. Conversely, every  $<$-component $\schm$ in $\forestf$ is of the type $T_\pathp$,  for any path $\pathp$ that intersects $\schm$.  
\end{corollary}

In the sequel we will be interested in the $<$-components of forests of the form $\tg{\tree}{\stem}$, for a given tree $\tree$.  By  $\tg{\schm}{\stem}_\ndu$ we will mean the $<$-component of $\tg{\tree}{\stem}$ that contains  the node $\ndu$.

\begin{lemma}\label{lem:schm=undiv}
Given a tree $\tree$ and a stem $\stem$  in it, let $\ndt, \ndu$ be nodes in $\tree$ such that $\stem < \ndt$ and $\stem < \ndu$. 
 Then $\tg{\schm}{\stem}_\ndt = \tg{\schm}{\stem}_\ndu$  if and only if $\patha \curlyvee_{\stem} \pathb$, for all paths $\patha$ and $\pathb$  containing respectively $\ndt$ and $\ndu$.
\end{lemma}

\begin{proof}
Let $\patha$ and $\pathb$ be paths containing $\ndt$ and $\ndu$, respectively. If $\tg{\schm}{\stem}_\ndt = \tg{\schm}{\stem}_\ndu$,  then (\ref{eq:schmerl-of-t}) from Proposition~\ref{prop:schmerl-comp} implies $\ndt \geqslant \ndt' \leqslant \ndu$ for some $\ndt'$ in $\tg{\trt}{\stem}$. Then $\ndt' \in\patha \cap \pathb $ and $\stem < \ndt'$, hence  $\patha \curlyvee_{\stem} \pathb$. 

Conversely, let $\patha \curlyvee_{\stem} \pathb$ and let $\ndt \in \patha$ and $\ndu \in \pathb$. 
Consider any $\ndv \in  \patha \cap \pathb$ with $\stem < \ndv$. 
We can assume w.l.o.g. that $\ndv \leqslant \ndt$ and $\ndv \leqslant \ndu$;  otherwise, $\ndt$ or $\ndu$ can play the role of $\ndv$. Then $\ndv \in \tg{\schm}{\stem}_\ndt \cap \tg{\schm}{\stem}_\ndu$, which implies $\tg{\schm}{\stem}_\ndt = \tg{\schm}{\stem}_\ndu$ by Proposition~\ref{prop:schmerl-comp}. 
\end{proof}

This observation shows that there is a one-to-one correspondence between equivalence classes modulo undividedness at $\stem$ and $<$-components of $\tg{\tree}{\stem}$. For every $[\patha]_{\curlyvee_{\stem}}$, the set $\{ \pathb \cap \tg{T}{\stem} : \pathb \in [\patha]_{\curlyvee_{\stem}} \}$ is the set of all paths of a unique $<$-component of $\tg{\tree}{\stem}$. Conversely, for every $<$-component $\schm$ of $\tg{\tree}{\stem}$, the set 
$\{ \pathp \cup  \stem : \pathp \ \text{is any path in} \ \schm \}$ 
is the set of all paths in a unique equivalence class modulo undividedness at $\stem$.

\begin{figure}[ht]
	\begin{center}
		\begin{picture}(130,140)
			\put(0,0){\line(1,1){130}}
			\qbezier(20,20)(40,40)(0,80)
			\qbezier(50,50)(70,70)(30,110)
			\qbezier(80,80)(100,100)(60,140)
			\put(0,130){{\large$\tree$}}
			\put(18,18){\circle*{3}}
			\put(35,45){\circle{2}}
			\put(40,50){\circle{2}}
			\put(45,55){\circle{2}}
			\put(65,75){\circle{2}}
			\put(70,80){\circle{2}}
			\put(75,85){\circle{2}}
			\put(95,105){\circle{2}}
			\put(100,110){\circle{2}}
			\put(105,115){\circle{2}}
			\put(124,113){\footnotesize{$\eta_0$}}
			\put(-5,65){\footnotesize{$\eta_{r_1}$}}
			\put(25,95){\footnotesize{$\eta_{r_2}$}}
			\put(55,125){\footnotesize{$\eta_{r_3}$}}
			\put(22,13){\tiny{$0$}}
		\end{picture}

		\caption{A tree with a stem that is branching$_2$ but not branching$_1$.\label{Fig:NonFinitelyBranchingTree}}
	\end{center}
\end{figure}
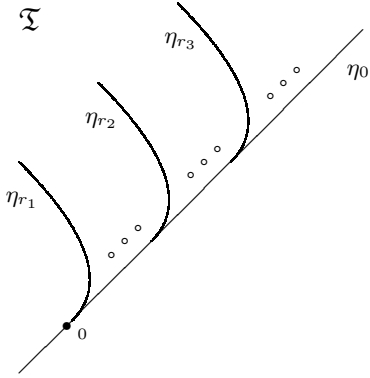

The undividedness classes at a given node represent the way in which a tree branches out at that node.  There are cases, though, in which this seems to conflict with the intuition. For instance, consider the tree $\tree$ obtained by taking a copy $\eta_0$ of the rationals and by attaching another copy $\eta_r$ of the rationals at every positive rational $r$ in $\eta_0$ (see Figure~\ref{Fig:NonFinitelyBranchingTree}).  
At the node 0 in $\eta_0$ of this tree there is only one undividedness class, but it can hardly be said that the tree does not branch out at 0. 

Then, two definitions can be considered.  
\begin{definition}\label{def: branch 1 and 2}
A stem $\stem$ in a tree $\trt$ is
\begin{enumerate}
\item  \defstyle{branching$_{1}$} iff  $\curlyvee_{\stem}$ has more than one equivalence class or, equivalently, iff  $\tg{\tree}{\stem}$ has more than one $<$-component or, equivalently, iff  there are two different paths $\patha_{1}, \patha_{2} \in \setpaths_{\stem}$ that are branching at $\stem$;

\item  \defstyle{branching$_{2}$} iff every node from $\tg{T}{\stem}$ has an incomparable node in $\tg{T}{\stem}$.  
\end{enumerate}
We say that a node in  $\trt$ is branching$_{i}$ iff the stem ending at that node is branching$_{i}$, for $i = 1,2$.
\end{definition}

\begin{proposition} 
\label{prop:branching} ~ 
\begin{enumerate}
\item {branching$_{1}$} implies {branching$_{2}$}. 
\item {branching$_{2}$} does not imply {branching$_{1}$}. 
\end{enumerate}
\end{proposition}

\begin{proof}
1. Let $\trt$ be any tree and $\stem$ be a  stem  in it. 
Suppose that $\stem$ is {branching$_{1}$} and let $\ndu$ be any node in $\tg{\trt}{\stem}$. Consider any $<$-component $\schm$ of  $\tg{\trt}{\stem}$ different from $\tg{\schm}{\stem}_\ndu$.   By Proposition~\ref{prop:schmerl-comp}, $\ndu \not\smile \ndv$ for all $\ndv \in \schm$.

2. The tree of Figure~\ref{Fig:NonFinitelyBranchingTree} provides a counterexample: the stem $(- \infty, 0]$ of  $\eta_0$ is branching$_{2}$ but not branching$_{1}$. 
\end{proof}

A \defstyle{bar} over a stem $\stem$ is a set $X$ of nodes such that $\stem < \ndv$ for every $\ndv \in  X$ and $X \cap \patha \not= \emptyset$ for every path $\patha$ containing $\stem$. A bar over the node $\ndt$ is a bar over the stem $ \tleq{T}{\ndt}$. Bars can also be used for describing how a tree branches out at a given stem, but the following result 
shows that nothing new is added in this way.

\begin{proposition}
A stem $\stem$ in the tree $\tree$ is branching$_2$ iff, for every bar $X$ over $\stem$, $|X| \geqslant 2$.
\end{proposition} 

\begin{proof}
Assume that $\stem$ is branching$_2$ and let $\ndv$ be any element of $\tg{T}{\stem}$. The set $\tg{T}{\stem}$ contains a node $\ndu$ that is not comparable with $\ndv$. Then, for every path $\pathp$ passing through $\ndu$, we have $\stem \subseteq \pathp$ and $\{\ndv\}  \cap \pathp = \emptyset$, and hence $\{\ndv\}$ cannot be a bar over $\stem$.

Conversely, suppose that there exists a node $t$ in $\tg{T}{\stem}$ such that $\ndt \smile \ndu$ for every $\ndu \in \tg{T}{\stem}$. Then, it can be easily verified that $\{ \ndt \}$ is a bar over 
$\stem$.
\end{proof}

\subsection{Bounded branching}

For any sets of nodes $X,Y$ in a tree, we say that \defstyle{$X$ underlies $Y$}, notation $X \leq_\undl Y$, if for every $\ndu \in Y$, there exists $\ndv \in X$ such that $\ndv \leqslant \ndu$.

\begin{definition}\label{def: n-branching 1 and 2}
For any $n \in \mathbb{N}$, the tree $\trt$ is
\begin{enumerate}
\item  \defstyle{$n$-branching$_{1}$ at the stem} $\stem$ iff  $\curlyvee_{\stem}$ has exactly $n$  equivalence classes or, equivalently, iff  $\tg{\tree}{\stem}$ has exactly $n$ $<$-components;
\item  \defstyle{$n$-branching$_{2}$ at the stem} $\stem$ iff  for every antichain $X$ in $\tg{\trt}{\stem}$, there exists an antichain $L_X$ in $\tg{\trt}{\stem}$ such that  $L_X \leq_\undl  X$ and $|L_X| = n$. 
\end{enumerate}

For $i \in  \{1,2\}$, $\tree$ is \defstyle{finitely branching$_i$} at $\stem$ if it is $n$-branching$_i$ at $\stem$ for some $n$.  If $\tree$ is not finitely branching$_i$ at $\stem$ then it is called \defstyle{infinitely branching$_i$} at $\stem$. 
\end{definition}

Recall that stems are assumed to be non-empty.  One can, however, expand the notion of $n$-branching$_{2}$ to include the `empty stem' as follows.  When $\stem$ is not a stem but the empty set then by definition $\tg{\trt}{\stem} = \tg{\trt}{\emptyset} = \trt$.  Hence we say that $\trt$ is \defstyle{$n$-branching$_{2}$ at $\emptyset$} when for every antichain $X$ in $\trt$ there exists an antichain $L_X$ in $\trt$ such that $L_X \leq_\undl X$ and $|L_X| = n$.  As before, $\trt$ is then called \defstyle{finitely branching$_2$} at $\emptyset$ when it is $n$-branching$_2$ at $\emptyset$ for some $n$.

\begin{definition}
	A tree $\trt$ is called \defstyle{finitely branching$_2$} when it is finitely branching$_2$ at $\emptyset$ and at each of its stems.
\end{definition}

It does not make sense to expand the of notion $n$-branching$_{1}$ to include the empty set, because the relation $\curlyvee_{\stem}$ is only defined for stems.  However, if one \textit{were} to attempt such an expanded definition, then from the fact that a tree $\trt$ is $n$-branching$_{1}$ at the stem $\stem$ if and only if $\tg{\tree}{\stem}$ has exactly $n$ $<$-components, a natural expansion would be that $\trt$ is $n$-branching$_{1}$ at $\emptyset$ if and only if $\trt$ has exactly $n$ $<$-components.  This would result in the trivial definition that each tree is $1$-branching$_{1}$ since trees, being connected, have exactly one $<$-component.

\begin{lemma}\label{lem:n-branching} Let $\trt$ be a tree and let $\stem$ be either the empty set or a stem in $\trt$.  If $\trt$ is finitely branching$_2$ at $\stem$, then every path in $\tg{\trt}{\stem}$  has an initial segment that is a bridge. 
\end{lemma} 
\begin{proof}

Let $\trt$ be $n$-branching$_2$ at $\stem$ and let $X$ be a maximal antichain in $\tg{T}{\stem}$.  Then, there exists a finite antichain $L \subseteq \tg{T}{\stem}$ such that $L \leq_\undl X$ and $|L| = n$. 
For every path $\pathp$ in $\tg{\trt}{\stem}$, let $\ndu_\pathp$ be the unique node from $L$ on $\pathp$. 

Now, fix any path $\patha$ in $\tg{\trt}{\stem}$.  
Then, for every $\ndw \in L$ we select a node $\ndw(\patha) \in \patha$ such that 
$\ndw(\patha)  \leqslant \ndw$ and $\ndw(\patha)  \leqslant \ndu_\patha$, if such node exists, else 
$\ndw(\patha) := \ndu_\patha$. 
Consider the set $L_\patha$ consisting of all these nodes. 
Let $\ndv_\patha$ be the least node in  $L_\patha$. It exists because $L_\patha$ is finite, linearly ordered, and $\ndu_\patha \in L_\patha$.   
Then the set $B_\patha = \{\ndv : \stem < \ndv \leqslant \ndv_\patha\}$  is an initial segment of $\patha$ and a bridge in $\patha$. Indeed, if there is a path $\pathp$  in $\tg{\trt}{\stem}$ such that 
$B_\patha \not\subseteq \pathp$ then $\pathp\neq \patha$ and 
$ \ndv_\patha \notin \pathp$, hence $ \ndu_\patha \notin \pathp$.
Suppose $B_\patha \cap \pathp \neq \emptyset$ and let $x \in B_\patha \cap \pathp$. 
Then $x < \ndu_\pathp$ and $x \leqslant \ndv_\patha \leqslant \ndu_\patha$.  
 By the choice of $\ndu_\pathp(A)$, it then must be the case that 
 $\ndu_\pathp(A)  \leqslant \ndu_\pathp$, hence 
 $\ndu_\pathp(A) \in \pathp$, hence $\ndu_\pathp(A) \leqslant \ndv_\patha$, hence $\ndu_\pathp(A) = \ndv_\patha$ by definition of $\ndv_\patha$. This contradicts the assumption that $B_\patha \not\subseteq \pathp$.
\end{proof}

The following results establish a relationship between the two notions of finite branching.

\begin{proposition}\label{prop:n-branching}
If the tree $\trt$ is $n$-branching$_2$ at the stem $\stem$, then $\tree$ is also $n$-branching$_1$ at $\stem$. 

\end{proposition}
\begin{proof}

For every $\patha$ in $\tg{\trt}{\stem}$, let $b_\patha$ be the intersection of all paths $\pathp$  in $\tg{\trt}{\stem}$ such that $\pathp \cap \patha \not= \emptyset$. 
It follows from Lemma~\ref{lem:n-branching} that each $b_\patha$ is a bridge. Assume $b_\patha \not= b_\pathb$ and $b_\patha \cap b_\pathb\not= \emptyset$, and consider $\ndu \in b_\patha \cap b_\pathb$ and $\ndv \in b_\patha \setminus b_\pathb$. Then the interval  $[\ndu, \ndv]$ is a furcation contained in $b_\patha$, which contradicts the fact that $b_\patha$ is a bridge.   Thus, for all paths $\patha$ and $\pathb$, either $b_\patha = b_\pathb$, or $b_\patha \cap b_\pathb = \emptyset$. Moreover, for $b_\patha \neq b_\pathb$, we have $\ndu \not\smile \ndv$ for all $\ndu \in b_\patha$ and $\ndv \in b_\pathb$.  Let $B$ be the set of all  bridges $b_\patha$. 

Let $X$ be any set contained in $\bigcup_{b \in B} b$ and such that $X \cap b$ is a singleton for each $b \in B$. Then $X$ in an antichain. Since $\tree$ is $n$-branching$_2$, there exists an antichain  $L_X \subseteq \tg{\trt}{\stem}$ such that $L_X \leqslant_\undl X$ and $|L_X| = n$. Then $L_X$ must contain exactly one element from each $b \in B$, hence $B$ consists of exactly $n$ bridges. 
Any two paths $\patha$ and $\pathb$  in $\setpaths_\stem$ (in $\tree$) are undivided at $\stem$ if and only if they contain the same bridge from  $B$, which is unique. Thus, there is a one-to-one correspondence between $B$  and the set of undividedness classes at $\stem$. Therefore $\tree$ is $n$-branching$_1$.
\end{proof}

The converse of this proposition does not hold. Indeed, the tree considered in Figure~\ref{Fig:NonFinitelyBranchingTree} has only one undividedness class at the node $0$ in $\eta_0$, but it is not finitely branching$_2$ at that node.  
The following result, which is a straightforward consequence of the previous proposition, establishes a more precise relationship between finite branching$_1$ and finite branching$_2$. 

\begin{corollary}
If the tree $\tree$ is $n$-branching$_1$ at the stem $\stem$, then either $\tree$ is $n$-branching$_2$ at $\stem$ or it is infinitely branching$_2$ at $\stem$.
\end{corollary}

Consider now, in a tree $\tree$, any path $\patha$ containing the stem $\stem$ and the undividedness class $[\patha]_{\curlyvee_{\stem}}$. Let $\tree'$ be the subtree of $\tree$, the domain $T'$ of which is the union of all paths in $[\patha]_{\curlyvee_{\stem}}$. The tree $\tree'$ is clearly $1$-branching$_1$ at $\stem$ and hence it is either $1$-branching$_2$ or infinitely branching$_2$ at $\stem$. 
This suggests that the only difference between the two notions derives from situations similar to that described in Figure~\ref{Fig:NonFinitelyBranchingTree}.

\section{Condensations of trees}
\label{sec:condensations}

The constructions of condensations of trees, introduced here, 
extend those of the condensations of linear orders as e.g. in \cite{Rosenstein}. 
 They relate to the notion of \emph{condensed trees}, formally defined further.

\subsection{Segments and bridges}
\label{subsec:bridges}

\begin{lemma}\label{lem:union_of_segms}
Let $\{\segm_i : i \in I \}$ be a set of segments in a given tree $\tree$, such that, for some  index $i_0$, $\segm_{i_0}\cap \segm_i \not= \emptyset$ for all $i \in I$. Assume also that $\segm^* = \bigcup_{i \in I} \segm_i$ is linearly ordered. Then $\segm^*$ is a segment.
\end{lemma}

\begin{proof}
We have to prove that $\segm^*$ is convex. Assume $x < z < y$, where $x$ and $y$ are elements of $\segm^*$. Assume  $x \in\segm_i$ and $y \in\segm_j$ and let $\ndt_i$ and $\ndt_j$ be elements of $\segm_i \cap\segm_{i_0}$ and $\segm_j \cap\segm_{i_0}$, respectively. Now several cases can be considered according the relative positions of $\ndt_i$, $\ndt_j$, $x$, $y$, and $z$ in $\segm_{i_0} \cup \segm_i \cup\segm_j$. In all cases it is easy to conclude $z \in \segm_{i_0} \cup \segm_i \cup\segm_j$ by using the convexity of $\segm_{i_0}$, or of $\in \segm_i$, or of $\segm_j$. Therefore, $z \in \segm^*$.
\end{proof}

\begin{proposition} \label{prop:union_of_bridges}
Let $\tree$ be a tree. Then: 
\begin{enumerate}
\item If $\{\segm_i : i \in I \}$ is a set of bridges such that  $\segm_{i_0}\cap \segm_i \not= \emptyset$ for some $i_0$ and all $i$ in $I$, then $\segm^*=\bigcup_{i \in I} \segm_i$ is a bridge.
\item 
	If $\segm$ is a bridge in $\tree$, then $\segm$ is contained in a unique maximal bridge.
\item 
The set of maximal bridges in $\tree$ forms a partition in $\tree$, i.e., the relation of
two nodes in $\tree$ belonging to the same maximal bridge forms an equivalence relation on $\tree$.
\end{enumerate}
\end{proposition}

\begin{proof}

	1. Let $\pathc$ be a path such that $\segm_{i_0} \subseteq \pathc$.  For every $i$, it holds that $\pathc \cap \segm_i \supseteq \segm_{i_0} \cap \segm_i \not= \emptyset$, which implies $\segm_i \subseteq \pathc$ because every $\segm_i$ is a bridge.  Then $\segm^* \subseteq \pathc $, which is linearly ordered, and hence, by Lemma \ref{lem:union_of_segms}, $\segm^*$ is a segment. 
	
Let $\pathp$ be any path such that  $\pathp \cap \segm^* \neq \emptyset$. Then, 
	$\pathp \cap \segm_{i} \neq \emptyset$ for some $i$.  Since $\segm_i$ is a bridge, we have $\pathp \supseteq \segm_{i} \supseteq \segm_{i} \cap \segm_{i_0}$, so that
	$\pathp \cap \segm_{i_0} \neq \emptyset$. Therefore, $\pathp \supseteq \segm_{i_0}$.  Then, by the assumption for $\segm_{i_0}$, 
	$\pathp$ has non-empty intersection with every $\segm_j$, and hence $\segm_j \subseteq \pathp$ for all $j$. Thus, $\segm^* \subseteq \pathp$.

		2.  Let $\mathcal{A}$ be the family of all bridges $\segmj$ such that $\segm \subseteq \segmj$. Then, by claim 1, $\segm^* =  \bigcup \mathcal{A}$ is a bridge. Clearly, it is the only maximal bridge containing $\segm$. 
	
	3. Let $\segm_1$ and $\segm_2$ be maximal bridges with non-empty intersection. Then, by claim 1, $\segm_1 \cup \segm_2$ is a bridge and hence, by the
	maximality of $\segm_1$ and $\segm_2$, we have $\segm_1 = \segm_1 \cup \, \segm_2 =\segm_2$.
\end{proof}

	For each $\ndt \in T$, the maximal bridge in $\tree = (T;<)$ containing $\ndt$ 
	(recall, that $\{ t \}$ is a bridge) 
	will be denoted as $\bridget$.	
Recall that the notation
$\bridget < \bridgeu$ 
means that $x < y$ for all $x \in \bridget$ and $y \in \bridgeu$, 
and $\bridget \smile \bridgeu$ will indicate that $\bridget < \bridgeu$ or $\bridget = \bridgeu$ or $\bridgeu < \bridget$. 

\begin{proposition} \label{ch2:Prop:OrderPreservedInCondensations}
	Let $\tree = (T;<)$ be a tree and let $\ndt,\ndu \in T$.
	\begin{enumerate}
		\item
			If $\ndt < \ndu$ and $\bridget \neq \bridgeu$ then $\bridget < \bridgeu$.
		\item
			If $\ndt \not\smile \ndu$ then $x \not\smile y$ for all $x \in \bridget$ and $y \in \bridgeu$.
	\end{enumerate}
\end{proposition}

\begin{proof}
	1. Let $\patha$ be a path in $\tree$ such that $\ndt, \ndu \in \patha$.  Since $\bridget$ and $\bridgeu$ are bridges then
	$\bridget,\bridgeu \subseteq \patha$ so that	$\bridget \cup \bridgeu$ is linearly ordered.  Since $\bridget \neq \bridgeu$ then
	$\bridget \cap \bridgeu = \emptyset$. Since bridges are segments, hence they are convex, it follows that $\bridget < \bridgeu$. 
	
	2. Follows from part 1.
\end{proof}

\begin{corollary} \label{ch2:Cor:OrderPreservedInCondensations}
	Let $\tree = (T;<)$ be a tree and let $\ndt, \ndu \in T$.  Then $\ndt \smile \ndu$ if and only if $\bridget \smile \bridgeu$.
\end{corollary}

\begin{proposition} \label{ch2:Prop:MaximalBridgeEquivalence}
	Let $\tree = (T;<)$ be a tree and $\ndt,\ndu \in T$.  The following are equivalent:
	\begin{enumerate}
		\item
			there exists a bridge $\segmb$ such that $\ndt,\ndu \in \segmb$;
		\item
			$\bridget = \bridgeu$;
		\item
			for every path $\pathp$ in $\tree$, $\ndt \in \pathp$ if and only if $\ndu \in \pathp$;
		\item
			for every node $\ndv \in \tree$, $\ndv \smile \ndt$ if and only if $\ndv \smile \ndu$.
	\end{enumerate}
\end{proposition}

\begin{proof}

	1 $\Longleftrightarrow$ 2: \ Immediate.

	2 $\Longrightarrow$ 3: \ Suppose $\bridget = \bridgeu$.  
		Let $\pathp$ be a path such that $\ndt \in \pathp$.  Then
	$\bridget \subseteq \pathp$, hence $\ndu \in \bridgeu\subseteq \pathp$.	
	Likewise, if $\ndu \in \pathp$ then $\ndt \in \pathp$.

	3 $\Longrightarrow$ 2: \ Suppose condition 3 holds.  If $\ndt = \ndu$ then the claim is
	immediate, so assume $\ndt \neq \ndu$ and let $\patha$ be a path with $\ndt \in \patha$.  Then $\ndu \in \patha$ and so $\ndt \smile \ndu$, say $\ndt < \ndu$.  
	Consider the segment $[\ndt,\ndu]$ and let $\pathp$ be any path such that 
	$\pathp \cap [\ndt,\ndu] \neq \emptyset$.  
	Then $\ndt \in \pathp$, so $\ndu \in \pathp$, hence 
	$[\ndt,\ndu] \subseteq \pathp$.  
	Therefore, $[\ndt,\ndu]$ is a bridge.  
	By Proposition \ref{prop:union_of_bridges}, $[\ndt,\ndu]$ is contained in a
	unique maximal bridge, hence $\bridget = \bridgeu$.

	3 $\Longrightarrow$ 4: \ Suppose condition 3 holds. 

	Let $\ndv \in \tree$ with $\ndv \smile \ndt$ and let $\patha$ be a path with $\ndv,\ndt \in \patha$.  Then $\ndu \in \patha$ and so $\ndv \smile \ndu$.  
	Likewise, if $\ndv \smile \ndu$ then $\ndv \smile \ndt$.

	4 $\Longrightarrow$ 3: \ Suppose condition 4 holds.   

		Let $\pathp$ be a path with $\ndt \in \pathp$.
	Since $\ndv \smile \ndt$ for every $\ndv \in \pathp$ then $\ndv \smile \ndu$ for every $\ndv \in \pathp$, hence $\ndu \in \pathp$.  
	Likewise, if $\ndu \in \pathp$ then $\ndt \in \pathp$.
\end{proof}

Thus, two nodes $x$ and $y$ in a tree belong to the same maximal bridge if and only if they satisfy in that tree the formula
\begin{equation} \label{ch2:Eqn:Beta(x,y)}
	\beta(x,y) := \forall z \left( z \smile x \leftrightarrow z \smile y \right).
\end{equation}
So, $\beta$ defines an equivalence relation on the set of nodes in the tree.

\subsection{Condensations} 
\label{subsec:Condensations-Ruaan}

\begin{definition}
\label{ch2:def:tree}
	Given a tree $\tree = (T;<)$, define the set $\bridgeT := \left\{ \bridget : t \in T \right\}$. The structure $\bridgetree :=  \left(\bridgeT;<\right)$ is called
	the \defstyle{condensation} (or, the  \defstyle{condensation quotient}) of the tree $\tree$.
\end{definition}

Thus, the condensation of a tree is its quotient structure generated by the equivalence relation of membership to the same maximal bridge. Consequently, the condensation of a tree shrinks all maximal bridges in it to single nodes.

The following is a straightforward consequence from  Proposition \ref{ch2:Prop:OrderPreservedInCondensations}.
\begin{proposition}
\label{ch2:prop:1}
	For any tree $\tree$, its condensation $\bridgetree$ is also a tree.
\end{proposition}

\begin{lemma} \label{ch2:Lem:MaximalBridgesAreSingletons}
	Let $\tree$ be a tree.  Then every bridge 
		in the tree $\bridgetree$ consists of a single node.
\end{lemma}

\begin{proof}
	Let $\bridget,\bridgeu \in \bridgetree$ and $\bridget \neq \bridgeu$.  Then $\ndt$ and $\ndu$ belong to different maximal
	bridges in $\tree$.  From Proposition \ref{ch2:Prop:MaximalBridgeEquivalence} we may conclude,
	without loss of generality, that there exists $\ndw \in \tree$ such that $\ndw \smile \ndt$ and
	$\ndw \not\smile \ndu$.  By Corollary \ref{ch2:Cor:OrderPreservedInCondensations} this implies that
	$\bridgew \smile \bridget$ and $\bridgew \not\smile \bridgeu$, so $\bridget$ and $\bridgeu$ belong to different
	maximal bridges in $\bridgetree$. Thus, all maximal bridges in $\bridgetree$ are singletons. 
\end{proof}

\medskip
\textbf{Remark 4.9}
\label{rem:condensation_via_Pt}
The condensation of a tree can be equivalently defined 
in terms of the sets $\setpaths_\ndt$  of paths containing a given node $\ndt$.
Let $\setpaths_T := \{\setpaths_\ndt : \ndt \in T\}$.  It is easily verified that the inverse inclusion $\supset$ is a tree relation on $\setpaths_T$. By 
Proposition \ref{ch2:Prop:MaximalBridgeEquivalence}, the nodes $\ndt$ and $\ndu$ belong to the same bridge if and only if $\setpaths_\ndt = \setpaths_\ndu$. Therefore, the map $\chi : \bridget \mapsto \setpaths_\ndt$ is a bijection. 
Moreover, $\bridget < \bridgeu$ if and only if $\setpaths_\ndt \supset \setpaths_\ndu$. Thus, the map $\chi : \bridget \mapsto \setpaths_\ndt$ is an isomorphism  from $\bridgetree$ to $(\setpaths_T; \supset)$.
\medskip

The operator $\bridge{\cdot}$ defines a canonical mapping
\begin{displaymath}
	\bridge{\cdot} : \, \tree \rightarrow \bridgetree.
\end{displaymath}
For $X \subseteq T$, $y \in \bridgeT$ 
and $Y \subseteq \bridgeT$, we denote
\begin{eqnarray*}
	\bridgeX & := & \left\{ \bridge{x} \in \bridgeT : x \in X \right\}, \\
	 \bridge{y}{^{-1}} & := & \left\{ x \in T : \bridge{x} = y \right\},  \\ 	
	 \bridgeY{^{^{-1}}} & := & \left\{ x \in T : \bridge{x} \in Y \right\}.
\end{eqnarray*}
Then $X \subseteq \bridge{\bridgeX}{^{^{^{\raisebox{.7ex}{\small -1}}}}}$ and $Y = \bridge{\bridgeY{^{^{-1}}}} $ for all $X \subseteq T$ and
$Y \subseteq \bridgeT$.

\begin{example} 
\label{ch2:ex:Condensation}
	Figure \ref{ch2:Fig:Condensation} shows a tree $\tree$ together with its condensation $\bridgetree$.  The
	bridges $\ndt_1$ through $\ndt_6$ are linear orders which may be infinite, and are condensed
	respectively to the nodes $\bridge{A_1}$ through $\bridge{A_6}$ in $\bridgetree$, so that $\bridgetree$ is finite.
	
	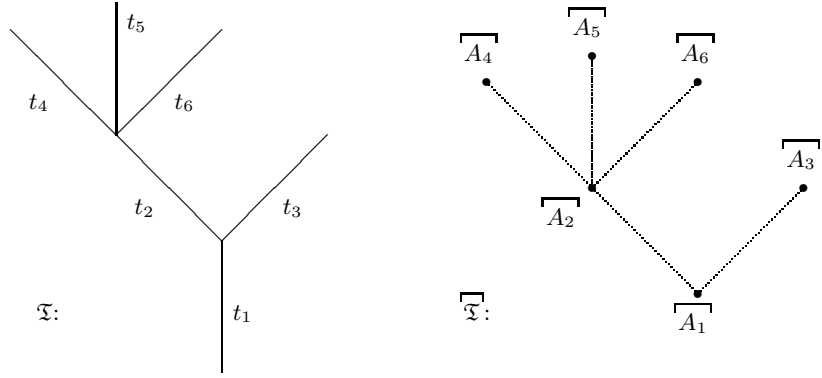
\begin{figure}[ht]
		\centering
		\begin{picture}(300,140)
			\put(80,-5){\line(0,1){50}}
			\put(80,45){\line(-1,1){40}}
			\put(80,45){\line(1,1){40}}
			\put(40,85){\line(-1,1){40}}
			\put(40,85){\line(0,1){50}}
			\put(40,85){\line(1,1){40}}
			\qbezier[40](260,25)(240,45)(220,65)
			\qbezier[40](260,25)(280,45)(300,65)
			\qbezier[40](220,65)(200,85)(180,105)
			\qbezier[50](220,65)(220,90)(220,115)
			\qbezier[40](220,65)(240,85)(260,105)
			\put(260,25){\circle*{3}}
			\put(220,65){\circle*{3}}
			\put(300,65){\circle*{3}}
			\put(180,105){\circle*{3}}
			\put(220,115){\circle*{3}}
			\put(260,105){\circle*{3}}
			\put(84,15){\footnotesize{$\ndt_1$}}
			\put(47,55){\footnotesize{$\ndt_2$}}
			\put(103,55){\footnotesize{$\ndt_3$}}
			\put(7,95){\footnotesize{$\ndt_4$}}
			\put(44,125){\footnotesize{$\ndt_5$}}
			\put(62,95){\footnotesize{$\ndt_6$}}
			\put(251,11){\footnotesize{$\bridge{A_1}$}}
			\put(201,51){\footnotesize{$\bridge{A_2}$}}
			\put(292,73){\footnotesize{$\bridge{A_3}$}}
			\put(170,113){\footnotesize{$\bridge{A_4}$}}
			\put(211,123){\footnotesize{$\bridge{A_5}$}}
			\put(252,113){\footnotesize{$\bridge{A_6}$}}
			\put(10,15){\footnotesize{$\tree$:}}
			\put(170,15){\footnotesize{$\bridgetree$:}}
		\end{picture}

		\caption{The condensation of a tree (see Example 
		\ref{ch2:ex:Condensation}). \label{ch2:Fig:Condensation}}
	\end{figure}
\end{example}

\begin{proposition} 
\label{ch2:Prop:StructureFromTreeToCondensation}
	Let $\tree = \treet$ be a tree and let $X \subseteq T$.
	\begin{enumerate}
		\item
			If $X$ is an antichain in $\tree$ then $\bridgeX$ is an antichain in $\bridgetree$. 
		\item
			If $X$ is linear in $\tree$  then $\bridgeX$ is linear in $\bridgetree$.
		\item
			If $X$ is a path 
			in $\tree$ then $\bridgeX$ is a path
			in $\bridgetree$.
		\item
			If $X$ is convex in $\tree$ then $\bridgeX$ is convex in $\bridgetree$.
		\item
			If $X$ is downward-closed in $\tree$ then $\bridgeX$ is downward-closed in $\bridgetree$. 
		\item
			If $X$ is upward-closed in $\bridgetree$ then $\bridgeX$ is upward-closed in $\bridgetree$.
	\end{enumerate}
\end{proposition}

\begin{proof}
	1., 2. From Corollary \ref{ch2:Cor:OrderPreservedInCondensations}.
	
	3. Let $X$ be maximal linear in $\tree$. Then $\bridgeX$ is linear.  Let
	$\bridgeu \in \bridgeT \backslash \bridgeX$.  Then $\ndu \not\in X$ so $\ndu \not\smile \ndv$ for some $\ndv \in X$.
	This implies $\bridgeu \not\smile \bridgev$ for $\bridgev \in  \bridgeX$. Thus, it follows that $\bridgeX$ is 	maximal linear in $\bridgetree$.

	4. Let $\bridget,\bridgeu \in \bridgeX$ and let $y \in \bridgeT$ be such that $\bridget < y < \bridgeu$.  Then
	$\ndt <  \bridge{y}{^{-1}} < \ndu$, hence  $ \bridge{y}{^{-1}}  \in X$.  This implies $ \bridge{\bridge{y}{^{^{-1}}}} = y \in \bridgeX$.

	The proofs of parts 5 and 6 are similar.
\end{proof}

\begin{proposition} \label{ch2:Prop:StructureFromCondensationToTree}
	Let $\tree = (T;<)$ be a tree and let $Y \subseteq \bridgeT$.
	\begin{enumerate}
		\item
			If $Y$ is linear  in $\bridgetree$ then $\bridgeY{^{^{-1}}}$ is linear  in $\tree$.
		\item
			If $Y$ is a path 
			 in $\bridgetree$ then $\bridgeY{^{^{-1}}}$ is a path 
			 in $\tree$.
		\item
			If $Y$ is convex in $\bridgetree$ then $\bridgeY{^{^{-1}}}$ is convex in $\tree$. 
		\item
			If $Y$ is downward-closed in $\bridgetree$ then $\bridgeY{^{^{-1}}}$ is downward-closed in $\tree$.
		\item
			If $Y$ is upward-closed in $\bridgetree$ then $\bridgeY{^{^{-1}}}$ is upward-closed in $\tree$.
	\end{enumerate}
\end{proposition}

	The proof is very similar to that of Proposition \ref{ch2:Prop:StructureFromTreeToCondensation}.

\medskip

Thus, paths, segments, stems, and branches are preserved between a tree and its condensation under the
mapping $\bridge{\cdot}$ and its inverse.

\subsection{Singular and emergent paths} 
\label{subsec:paths}
 
\begin{definition}
	A path $\patha$ in a tree $\tree$ is called	\defstyle{singular} 
	if there exists $\ndt \in \patha$ such that $\tgeq{T}{\ndt}$ is
	linear (and, therefore, $\tgeq{T}{\ndt} \subseteq \patha$).  Otherwise $\patha$ is called \defstyle{emergent}.
\end{definition}
Singular paths are of special interest from a model-theoretic viewpoint because each singular path can be defined by a first-order formula: if $\patha$ is singular and $\ndt \in \patha$ is such that $\tgeq{T}{\ndt}$ is linear then $\patha$ can be defined in $\trt$ by the formula $\varphi(x)$ given by $x \smile \ndt$. 

\begin{proposition} 
\label{ch2:Prop:SingularWhenTerminal}
	A path $\patha$ in a tree $\tree$ is singular if and only if 
	$\bridge{\patha}$
	contains a greatest node in $\bridgetree$. 
\end{proposition}

\begin{proof}
	Let $\patha$ be singular and let $\ndt \in \patha$ be such that $\tgeq{T}{\ndt}$ is linear. Then note that $\tgeq{T}{\ndt}$ is a bridge in $\tree$ and  $\tgeq{T}{\ndt} \subseteq \bridget$, hence $\bridget$ is the greatest node of $\bridge{\patha}$

	Conversely suppose $\bridge{\patha}$ contains a greatest node $\bridget$.  
Then $\ndt \in \patha$.  
Let $\ndu,\ndv \in \tgeq{T}{\ndt}$.  
 Then $\bridget \leqslant \bridgev$ and 
 $\bridget \leqslant \bridgev$, hence 
$\bridgeu = \bridgev = \bridget$, so 
	$\ndu \smile \ndv$. Thus, $\tgeq{T}{\ndt}$ is linear.
\end{proof}

Thus, a path $\patha$ is emergent if and only if $\bridge{\patha}$ does not contain a greatest node.

\smallskip
A tree $\tree$ is called \defstyle{well-founded} 
when every non-empty set of nodes from $\tree$ contains a minimal node.
Note that a tree is well-founded if and only if each of its paths is well-ordered. 

\begin{proposition}
\label{Prop:well-founded}
	Let $\tree$ be a well-founded tree and let $\patha$ be a path in $\tree$. Then: 
		\begin{enumerate}
		\item
		$\bridgetree$ is a well-founded tree;
		\item
			$\patha$ is singular if and only if the order type of $\bridge{\patha}$ is a successor ordinal;
		\item
			$\patha$ is emergent if and only if the order type of $\bridge{\patha}$ is a limit ordinal.
	\end{enumerate}
\end{proposition}

\begin{proof}
	First, every path in $\bridgetree$ is a quotient $\bridge{\pathp}$ of a path $\pathp$ in $\tree$, and $\bridge{\pathp}$ is well-ordered (hence, has the order type of an ordinal) whenever $\pathp$ is well-ordered. The claims now follow from Proposition \ref{ch2:Prop:SingularWhenTerminal}.
\end{proof}

\subsection{Condensed trees}

\begin{definition}
	A tree $\tree$ is called \defstyle{condensed} when $\tree \cong \bridgetree$.
\end{definition}

\begin{proposition} \label{ch2:Prop:CondensedTreeEquivalences}
	Let $\tree = (T;<)$ be a tree.  The following conditions are equivalent:
	\begin{enumerate}
		\item
			$\tree$ is condensed;
		\item
			$\tree \cong \bridge{\tree'}$ for some tree $\tree'$;
		\item
			$\bridgew = \{ \ndw \}$ for every $\ndw \in T$;
		\item
			$\setpaths_\ndu \not= \setpaths_\ndv$ for all $\ndu \not= \ndv$ in $T$.\footnote{Because of this property, condensed trees are called \defstyle{totally branching} in \cite{SabbadinZanardo}.}
	\end{enumerate}
\end{proposition}

\begin{proof}
	1 $\Longrightarrow$ 2: \ Let $\tree$ be condensed.  Then $\tree \cong \bridgetree$.

	2 $\Longrightarrow$ 3: \
	Suppose that $f : \tree \to \bridge{\tree'}$ is an isomorphism.  Let $\ndt,\ndu \in T$ with $\ndt \neq \ndu$.  Then
	$f(\ndt) \neq f(\ndu)$ and, by Lemma \ref{ch2:Lem:MaximalBridgesAreSingletons}, $f(\ndt)$ and $f(\ndu)$ belong to different maximal bridges in $\bridge{\tree'}$.  From Proposition~\ref{ch2:Prop:MaximalBridgeEquivalence} we may conclude, without loss of generality, the existence of a node $\nds$ in $\bridge{\tree'}$ such that $\nds \smile f(\ndt)$ and $s \not\smile f(\ndu)$.  Since $f$ is an isomorphism then $f^{-1}(\nds) \smile \ndt$ and $f^{-1}(\nds) \not\smile \ndu$.  Hence by Proposition~\ref{ch2:Prop:MaximalBridgeEquivalence}, $\bridget \neq \bridgeu$ and the result follows.

	3 $\Longrightarrow$ 1: \ If $\bridgew = \{ \ndw \}$ for every $\ndw \in T$, then the canonical map given as $\ndw \mapsto \bridgew$ defines an isomorphism from $\tree$ to $\bridgetree$.

	3 $\Longleftrightarrow$ 4: \
	In Remark 4.9 
	we have observed that the map $\bridget \mapsto \setpaths_\ndt$ is a bijection. Then every bridge is a singleton if and only if part 4 holds.
\end{proof}

\begin{proposition}\label{prop:non-leaf nodes}
Let $\tree$ be a tree and let $y$ be any non-leaf node in $\bridgetree$. Then $\bridge{y}{^{-1}}$ is a final segment of a branching$_2$ stem $\stem$ in $\tree$, i.e. a subset of $\stem$ which is upward-closed in $\stem$. 
\end{proposition}

\begin{proof} 
The set $\bridge{y}{^{-1}}$ is a segment. The set $\stem = \{\ndt: \ndt \leqslant \ndu \textrm{ for some } \ndu \in \bridge{y}{^{-1}} \}$ is a stem in $\tree$, and $\bridge{y}{^{-1}}$ is a final segment of it. Consider any node $\ndv$ in $\tg{T}{\stem}$. If $\ndv \smile \ndt$ for every $\ndt \in \tg{T}{\stem}$, then  $\bridge{y}{^{-1}} \cup \{\ndu : \bridge{y}{^{-1}} < \ndu \leqslant \ndv \}$ is a bridge, which contradicts the maximality of $\bridge{y}{^{-1}}$.
Thus, $\stem$ is a branching$_2$ stem. 
\end{proof}

\begin{proposition}\label{prop:fin branch -> well founded}
Let $\tree$ be a finitely branching$_2$ tree. Then $\bridgetree$ is well-founded. 
\end{proposition}

\begin{proof} Suppose $\bridgetree$ is not well-founded. Then there exists an infinite  sequence $\{y_i\}_{i\in \mathbb{N}} \subseteq \bridgeT$ such that $y_{i+1} < y_i$ for all $i\in \mathbb{N}$. 

Let $Y = \{y \in \bridgeT : y < y_i \mbox{ for all } i\in \mathbb{N} \}$ (note that $Y$ may be empty)
and $Z = \{y \in \bridgeT : Y < y < y_0\}$. 
Then either $\bridgeY{^{^{-1}}} = \emptyset$ or $\bridgeY{^{^{-1}}}$ is a stem in $\tree$, and $\bridge{Z}{^{^{-1}}}$ is a segment in $\tree$ such that there is no node $\ndv \in T$ for which $\bridgeY{^{^{-1}}} < \ndv < \bridge{Z}{^{^{-1}}}$. Then every initial subsegment 
of  $\bridge{Z}{^{^{-1}}}$ 
(i.e., a downward-closed subset of  $\bridge{Z}{^{^{-1}}}$) 
contains infinitely many disjoint bridges and hence it is a furcation. 
From Lemma~\ref{lem:n-branching} it follows that $\bridgetree$ is not finitely branching$_2$ from  $\bridgeY{^{^{-1}}}$. 
\end{proof}

Some useful observations follow from the propositions above: 

\begin{enumerate}
\item A tree is condensed iff each of its (maximal) bridges is a singleton. 

\item A tree is condensed iff each of its non-leaf nodes is a branching$_2$ node. 

\item The condensation of a condensed tree is (isomorphic to) the same tree. 
\end{enumerate}

\subsection{Refinements and homeomorphisms of trees}
\label{subsec:Subdivisions} 

 A \defstyle{refinement} 
of a tree $\tree$ is again a tree obtained from $\tree$ by inserting nodes inside bridges. (Note the class of trees is closed under such insertions.) 
Two trees are  \defstyle{homeomorphic} if they  have isomorphic refinements.  The relation between trees of being homeomorphic is an equivalence relation.  This fact is readily seen using the following observation.  

Trees $\trt$ and $\trs$ are homeomorphic if and only if $\bridge{\trt} \cong \bridge{\trs}$.  The forward direction is clear: if $\trt'$ and $\trs'$ are isomorphic refinements of $\trt$ and $\trs$ then $\bridge{\trt} \cong \bridge{\trt'} \cong \bridge{\trs'} \cong \bridge{\trs}$.  Conversely, if $\bridge{\trt} \cong \bridge{\trs}$ then a tree $\trt'$ that is isomorphic to refinements of both $\trt$ and $\trs$ can be constructed as follows.  Let $f : \bridge{\trt} \to \bridge{\trs}$ be an isomorphism.  Obtain $\trt'$ from $\trt$ as follows.  For each $\bridge{\ndt}$ in $\bridge{\trt}$, if $\bridge{\ndt} \cong f\left(\bridge{\ndt}\right)$ then leave the bridge $\bridge{\ndt}$ in $\trt$ as it is, while if $\bridge{\ndt} \not\cong f\left(\bridge{\ndt}\right)$ then replace the bridge $\bridge{\ndt}$ in $\trt$ with a copy of $\bridge{\ndt}$ followed by a copy of $f\left(\bridge{\ndt}\right)$.  It is readily seen that this tree $\trt'$ is isomorphic to refinements of both $\trt$ and $\trs$.

Now, take a tree $\tree$ and consider a partition of it into a set $S$ of bridges. The ordering in $\tree$ is inherited in $S$, just like in its condensation quotient, which turns the set $S$ into a tree itself, denoted by $\tree(S)$.  
 The canonical mapping from $\tree$ to  $\tree(S)$ is called a \defstyle{homeomorphism of  $\tree$ onto $\tree(S)$}, and $\tree(S)$ is called a \defstyle{homeomorphic abstraction} of the tree $\tree$.  
 Thus, two trees are homeomorphic if and only if they share a common, up to isomorphism, homeomorphic abstraction. Note that homeomorphic abstractions generalise the condensation quotient construction, obtained when $S$ is the set of maximal bridges, and many of the results about condensations in this section apply likewise to homeomorphic abstractions. We leave out the routine details.  

Clearly, every homeomorphism of trees is a homomorphism of partial orders, but not vice versa.  
Indeed, any tree-like partial order can be homomorphically mapped onto a linear order, but that, in general, is not a homeomorphism of trees. 

We also remark that the notion of homeomorphism of trees has a topological nature, in the following sense. The set $\setpaths(\trt)$ of paths in a tree $\trt = (T;<)$ can be endowed 
with a topological structure by considering the family $\setpaths_T$ of all sets of the type 
$\setpaths_\ndt$, for $t \in T$,  
as a subbase of open sets. Denote this topology by $\tau_\trt$. It turns out that $(\setpaths(\trt), \tau_\trt)$ is a {\sl non-Archimedean} topological space (\cite{Nyikos}) because  it is Hausdorff and  $\{\setpaths_\ndt : \ndt \in T\}$ is a {\em rank~1} subbase, that is, for all $\ndt, \ndt' \in T$, either $\setpaths_\ndt \subseteq \setpaths_{\ndt'}$, or $\setpaths_{\ndt'} \subseteq \setpaths_{\ndt}$, or $\setpaths_\ndt \cap \setpaths_{\ndt'} = \emptyset$. Thus, by Remark \ref{rem:condensation_via_Pt}, it follows that homeomorphic trees generate homeomorphic topologies, which also justifies the adopted terminology.

\section{Condensed forkings and extensions} 
\label{sec:Condensed extensions} 
The definition of condensed trees can be extended to \emph{condensed forests} in a natural way.   
It can be easily verified that all results for condensed trees in the previous section hold for condensed forests, too. Note also that a forest is condensed if and only if all of its $<$-components are condensed trees.

Propositions~\ref{ch2:Prop:StructureFromTreeToCondensation} and \ref{ch2:Prop:StructureFromCondensationToTree} show that the condensation of a given tree $\tree$ produces a new tree having the same `branching structure' as $\tree$. But most of the inner structure of its paths is lost with this operation. It is sometimes desirable to preserve the type and structure of the paths, while ensuring condensed branching. For that purpose, 
following \cite{SabbadinZanardo}, here we introduce and study an alternative construction which does not condense the tree, but, instead, extends it to produce a new condensed tree in which every path is isomorphic to a path in the original tree. 

\smallskip
We will present two versions of the construction, a full version and a refined one.
The basic idea of both 
is that non-singleton bridges can be eliminated by duplicating nodes and the subtrees generated by those nodes.  Furthermore, the property of being condensed is ensured by duplicating the subtrees rooted \emph{at every node of the original tree} in the full version, and \emph{only at the non-branching nodes of the original tree} in the refined version.

As the construction (in both versions) generally produces forests, we will present it applied not just to trees, but to forests. Every connected $<$-component of each of the resulting extensions will be a  condensed tree.  

We first define the full version.  
Given any $\ndt \in T$, we start duplicating all nodes in $\tleq{T}{\ndt}$. This is formally obtained by considering functions that assign $0$ or $1$ to each element of $\tleq{T}{\ndt}$. The new structure will consist of all these functions for $\ndt$ ranging over $T$.  
We define
\begin{equation}\label{ch?:eq:CondensedExt} \conextset := \bigcup \left\{\,2^{{\tleq{T}{\ndt}}} :\ndt \in T \right\}   
\end{equation}
where $2^X$ denotes the set of all functions $f: X \to \{0,1\}$.  
When $X \subset Y$, $f \in 2^X$, and $g \in 2^Y$, by $f \subset g$ we denote the claim that $g$ extends $f$, i.e., that $f$ is the restriction $g{\upharpoonright}_X$ of $g$ to $X$. 

For any given $\tree$, $\conext$ will denote the structure $(\conextset;\subset)$, defined above. 

\smallskip
Now, the definition of $\conext$ can be refined in order to obtain a smaller extension 
$\conextset^\dagger$ according to the following informal explanation. Observe first that the duplication of nodes in the construction of $\conextset$ essentially  aims at turning non-singleton bridges into furcations. Moreover, if $\ndt$ is the initial node of a bridge, then, for every $\ndu < \ndt$, the interval $[\ndu, \ndt]$ is a furcation, and hence the duplication of $\ndt$ produces only copies of already existing furcations. This means that the duplication of $\ndt$ is unnecessary. Then, we set
\begin{equation}\label{ch?:eq:CondensedExt-} \conextset^\dagger := \conextset \setminus \{ f  : f(\ndt) = 0 \textrm{ for some initial node } \ndt \textrm{ of a maximal bridge}  \}.
\end{equation}
Thus, $\conextset^\dagger$ is the set of all functions in $\conextset$ such that $f(\ndt) = 1$ whenever $\ndt$ is the initial point of a maximal bridge. 
We now define $\conext^\dagger := (\conextset^\dagger ; \subset)$. 

Observe that, if $\tree$ is condensed, the set $\{\ndt\}$ is a maximal bridge for every node $\ndt$, which is its initial node. Then, $\conextset^\dagger$ is the set of all functions in $\conextset$ that constantly take the value $1$. For every node $\ndt$, there is only one function of this kind having  $\tleq{T}{\ndt}$ as domain and  hence $\conext^\dagger$ is (a tree) isomorphic to $\tree$. 
Also, if $\tree$ is not rooted, then $\conext^\dagger$ might be a (condensed) proper forest.

\begin{example} 
\label{ex2:CondensedExtension} 
Figure~\ref{fig:ex2:CondensedExtension} shows a tree $\tree$ on the left, 
a $<$-component of the full construction $\conext$ applied to $\tree$ in the middle, 
and the refined version $\conext^\dagger$ on the right.  

	\begin{figure}[ht]
		\centering
\begin{picture}(330,80)
			\qbezier[40](30,5)(30,25)(30,35)
			\put(30,5){\circle*{3}}
			\put(30,35){\circle*{3}}
			\qbezier[40](30,35)(20,50)(10,65)
			\put(10,65){\circle*{3}}
			\qbezier[40](30,35)(40,50)(50,65)
			\put(50,65){\circle*{3}}

			\put(34,3){\footnotesize{$\ndt$}}
			\put(34,30){\footnotesize{$\ndu$}}
			\put(9,70){\footnotesize{$\ndv$}}
			\put(43,70){\footnotesize{$\ndw$}}
			
			\put(0,15){\footnotesize{$\tree$:}}
			
				\qbezier[40](150,5)(135,20)(120,35)
			\put(150,5){\circle*{3}}
                   \put(120,35){\circle*{3}}
                   
			\qbezier[40](120,35)(130,50)(140,65)
			\put(140,65){\circle*{3}}
			\qbezier[40](120,35)(110,50)(100,65)
			\put(100,65){\circle*{3}}
			
			\qbezier[40](150,5)(165,20)(180,35)
			\put(180,35){\circle*{3}}
			\qbezier[40](180,35)(170,50)(160,65)
			\put(160,65){\circle*{3}}
			\qbezier[40](180,35)(190,50)(200,65)
			\put(200,65){\circle*{3}}
	
			\put(154,3){\footnotesize{$\ndt_0$}} 
			\put(184,30){\footnotesize{$\ndu_{01}$}}
			\put(125,30){\footnotesize{$\ndu_{00}$}}
			
			\put(70,70){\footnotesize{$\ndv_{000}$}}
			\put(108,70){\footnotesize{$\ndw_{000}$}}
			\put(89,70){\footnotesize{$\ndv_{001}$}}
			\put(130,70){\footnotesize{$\ndw_{001}$}}
			\put(150,70){\footnotesize{$\ndv_{010}$}}
			\put(188,70){\footnotesize{$\ndw_{010}$}}
			\put(169,70){\footnotesize{$\ndv_{011}$}}
			\put(210,70){\footnotesize{$\ndw_{011}$}}
			
			\put(80,15){\footnotesize{$\forestf_{\tree_0}$:}}
			
			\qbezier[40](120,35)(100,50)(80,65)
			\put(80,65){\circle*{3}}
			\qbezier[40](120,35)(120,50)(120,65)
			\put(120,65){\circle*{3}}
			
			\qbezier[40](180,35)(180,50)(180,65)
			\put(180,65){\circle*{3}}
			\qbezier[40](180,35)(200,50)(220,65)
			\put(220,65){\circle*{3}}

				\qbezier[40](300,5)(285,20)(270,35)
			\put(300,5){\circle*{3}}
                   \put(270,35){\circle*{3}}
                   
			\qbezier[40](270,35)(280,50)(290,65)
			\put(290,65){\circle*{3}}
			\qbezier[40](270,35)(260,50)(250,65)
			\put(250,65){\circle*{3}}
			
			\qbezier[40](300,5)(315,20)(330,35)
			\put(330,35){\circle*{3}}
			\qbezier[40](330,35)(320,50)(310,65)
			\put(310,65){\circle*{3}}
			\qbezier[40](330,35)(340,50)(350,65)
			\put(350,65){\circle*{3}} 
	
			\put(304,3){\footnotesize{$\ndt$}} 
			\put(334,30){\footnotesize{$\ndu_1$}}
			\put(275,30){\footnotesize{$\ndu_0$}}
			\put(250,70){\footnotesize{$\ndv_0$}}
			\put(310,70){\footnotesize{$\ndv_1$}}
			\put(278,70){\footnotesize{$\ndw_0$}}
			\put(338,70){\footnotesize{$\ndw_1$}}
			\put(240,15){\footnotesize{$\conext^\dagger$:}}
	
		\end{picture}

		\caption{A tree $\tree$, a $<$-components of $\conext$, and $\conext^\dagger$. \label{fig:ex2:CondensedExtension}}
	\end{figure}
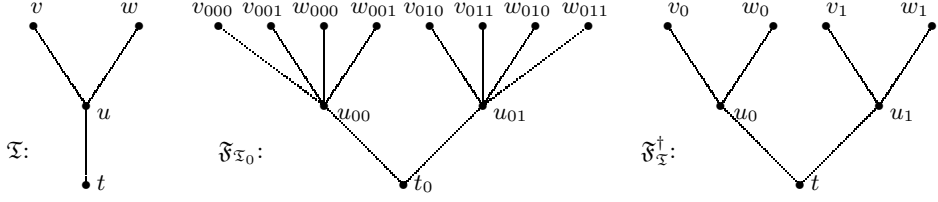
\end{example}

As another example, if $\tree = \left(\mathbb{N};<\right)$, then $\conext^\dagger$ is the infinite binary branching tree and $\conext$ consists of two copies of that tree.

\begin{proposition}[\cite{SabbadinZanardo}, Proposition~5.19]~\label{ch4:condensed forests} \\
	For any tree $\trt = \treet$, the structure $\conext = (\conextset;\subset)$ is a condensed forest.
\end{proposition}

The same holds for the refined extension $\conext^\dagger$.

\begin{proposition}\label{ch?:prop:CondensedExt-}
	For any tree $\trt = \treet$, the structure $\conext^\dagger = (\conextset^\dagger; \subset)$ is a condensed forest.
\end{proposition}

\begin{proof}
The definition of forest does not involve existential assumptions. Then $\conext^\dagger$  is a forest because $\conextset^\dagger \subseteq \conextset$ and $\conext$ is a forest. We show that  $\setpaths_f \not= \setpaths_g$, for all different nodes $f$ and $g$ in $\conext^\dagger$. This implies that $\conext^\dagger$ is condensed by Proposition~\ref{ch2:Prop:CondensedTreeEquivalences}.

The claim is trivial if $f \not\smile g$. Assume $f \subset g$ and let $\ndt$ be an element of $\dom{g} \setminus \dom{f}$. Two cases can be considered.

\smallskip
Case a: $\ndt$ is not the initial node of a maximal bridge. Let $g'$ be the element of $\conextset^\dagger$ defined by: $\dom{g'} = \dom{g}$; $g'(\ndu) = g(\ndu)$ for all $\ndu \not= \ndt$ in $\dom{g}$; $g'(\ndt) \not= g(\ndt)$. Then $g$ and $g'$ are $\subset$-incomparable nodes and $f \subset g'$. Any path containing $g'$ in  $\conext^\dagger$ is an element of $\setpaths_f \setminus \setpaths_g$.

\smallskip
Case b: $\ndt$ is the initial node of a maximal bridge $b$. Assume $f \in 2^{\tleq{T}{\ndu}}$, so that $\ndu < \ndt$. Observe that, if $\ndv \smile \ndt$  for all $\ndv > \ndu$, then $b \, \cup [\ndu, \ndt]$ is  a bridge and this contradicts the maximality of $b$. Then we can consider a node $\ndv > \ndu$ such that $\ndv \not\smile \ndt$. Let $g'$ be any extension of $f$ in $2^{\tleq{T}{\ndv}}$. Every path containing $g'$ in  $\conext^\dagger$ is an element of $\setpaths_f \setminus \setpaths_g$.
\end{proof}

Now, let $\fftt$ denote the function from  $\conextset$ onto $T$ defined by 
\begin{equation}\label{ch4:function from F_T to T}
\fftt(f) = \ndt \ \ \textrm{for all } f \in  2^{\tleq{T}{\ndt}}.
\end{equation}

\begin{proposition}\label{ch?:prop:CondensedExt} For any tree $\trt = \treet$, the function $\fftt$ has the following properties:
\begin{enumerate}
\item $\fftt$ is order-preserving;

\item for every path $\tilde\pathp$ in $\conext$, $\fftt\left(\tilde\pathp\right)$ is a path in $\tree$ and the restriction of $\fftt$ to $\tilde\pathp$ is an isomorphism;

\item every path in $\tree$ is the $\fftt$-image of a path in $\conext$.
\end{enumerate}
\end{proposition}

\begin{proof} 
Observe first that, for $f \in2^{{\tleq{T}{\ndt}}}$ and $g  \in 2^{{\tleq{T}{\ndu}}}$, $f \subset g$ iff $\ndt < \ndu$ and $f$ is the restriction $g{\upharpoonright}_{\tleq{T}{\ndt}}$ of $g$ to $\tleq{T}{\ndt}$. Then part 1 holds.

Consider a path $\tilde\pathp$ in $\conext$ and let $\patha$ be $\fftt(\tilde\pathp)$. The linearity of $\subset$ on $\tilde\pathp$ implies that $\patha$ is linearly ordered by $<$. Moreover, $\bigcup \{ f : f \in \tilde\pathp \}$ is a function $\xi$ in $2^\patha$ and every element of $\tilde\pathp$ can be written as $\xi{\upharpoonright}_{\tleq{T}{\ndu}}$ for some $\ndu \in \patha$. Assume that $\patha \cup \{\ndv\}$ is linearly ordered.
If $\ndv < \ndu$ for some $\ndu \in \patha$ then $\tilde\pathp \cup \left\{\xi{\upharpoonright}_{\tleq{T}{\ndv}}\right\}$ is also linearly ordered. So $\xi{\upharpoonright}_{\tleq{T}{\ndv}} \in \tilde\pathp$ and $\ndv \in \patha$. If $\ndu < \ndv$ for all $\ndu \in \patha$, so that $\ndv \not\in \patha$, consider the function $\xi' = \xi \cup \{ \langle \ndv, 0 \rangle\}$. For every $f \in \tilde\pathp$, we have $f \subset \xi'{\upharpoonright}_{\tleq{T}{\ndv}}$ and hence the set $\tilde\pathp \cup  \left\{\xi'{\upharpoonright}_{\tleq{T}{\ndv}}\right\}$ is linearly ordered by $\subset$. The maximality of $\tilde\pathp$ implies that it contains $\xi'{\upharpoonright}_{\tleq{T}{\ndv}}$, which contradicts $\ndv \not\in \patha$. Then  $\patha$ is a path in $\tree$.

We have observed that  $\tilde\pathp = \{ \xi{\upharpoonright}_{\tleq{T}{\ndt}} : \ndt \in \patha\}$. This implies that the restriction of $\fftt$ to $\tilde\pathp$ is injective and surjective on $\patha$. Then it is an isomorphism because $\fftt$ is order preserving.  This concludes the proof of part 2.

\smallskip
Given any path $\pathp$ in $\tree$, consider an element $\chi$ of $2^\pathp$ and set 
\begin{equation}\label{ch?:eq:CondensedExt1}
\chi^* := \{ \chi{\upharpoonright}_{\tleq{T}{\ndt}} : \ndt \in \pathp\}.
\end{equation}
Observe first that $\pathp = \fftt(\chi^*)$.  It is easily verified that $\chi^*$ is linearly ordered by $\subset$. If $f \in \, 2^{\tleq{T}{\ndt}}$ for some $\ndt \in T$ and $\chi^* \cup \{f\}$ is linearly ordered by $\subset$, then the maximality of $\pathp$ implies $ \ndt \in \pathp$. Since $f \smile \chi{\upharpoonright}_{\tleq{T}{\ndt}}$ and $f \in \, 2^{\tleq{T}{\ndt}}$, then $f = \chi{\upharpoonright}_{\tleq{T}{\ndt}} \in \chi^*$. We conclude that $\chi^*$ is a path in $\conext$ and part 3 holds.
\end{proof}

Note that $\conext$ is a proper forest for every tree $\tree$. Indeed, consider any node $\ndt$ in $\tree$ and let $f$ and $g$ be elements of $2^{\tleq{T}{\ndt}}$ such that, for every $\ndu \in \tleq{T}{\ndt}$, $f(\ndv) \not= g(\ndv)$ for some $\ndv \leqslant \ndu$. For no $h \in \conextset$ do we have $h \subset f$ and $h \subset g$. 

In the particular case where $\tree$ is rooted, $\conext$ consists of two disjoint condensed trees. In fact, if $\ndt_0$ is the root, for every $f \in 2^{\tleq{T}{\ndt}}$, either $\{\langle \ndt_0, 0 \rangle\} \subseteq f$ or $\{\langle \ndt_0, 1 \rangle\}\subseteq f$. Then, $\{\langle \ndt_0, 0 \rangle\}$ and $\{\langle \ndt_0, 1 \rangle\}$ are the roots of the two disjoint subtrees, $\tree_0$ and $\tree_1$.
Observe that $T_0 = \{f \in \conextset : \{\langle \ndt_0, 0 \rangle\} \subseteq f\}$ and $T_1 = \{f \in \conextset : \{\langle \ndt_0, 1 \rangle\} \subseteq f\}$. Then $\tree_0$ and $\tree_1$ are isomorphic. 

These observations can be generalized as follows.

\begin{lemma}\label{ch4: paths in F_T}
For every path $\pathp$ in the tree $\tree$  and every $<$-component $\tree^*$ of $\conext$, there exists a path $\tilde\pathp$ in  $\tree^*$ such that $\fftt(\tilde\pathp) = \pathp$. 

Likewise for every $<$-component of $\conext^\dagger$. 
\end{lemma}

\begin{proof}
We will only prove the claim for $\conext$, as the proof for $\conext^\dagger$ is very similar. 
Consider any path $\tilde\patha$ in  $\tree^*$ and let $\patha$ be $\fftt\big(\tilde\patha\big)$. Let $\stem = \pathp \cap \patha$. Since the restriction of $\fftt$ to $\tilde\patha$ is an isomorphism, $\stem$ is the $\fftt$-image of a downward-closed subset $\tilde\stem$ of $\tilde\patha$ in $\tree^*$. The set $\xi = \cup \{ f : f \in \tilde\stem \}$ is a function from $\stem$ to $2$ and it can be extended to a function $\chi$ from $\pathp$ to 2. The set $\{\chi{\upharpoonright}_{\tleq{T}{\ndt}} : \ndt \in \pathp\}$ is a path $\tilde\pathp$ in $\tree^*$ and $\fftt\big(\tilde\pathp\big) = \pathp$.
\end{proof}

\begin{proposition}\label{ch4: isomorphic components}
All $<$-components of $\conext$ are isomorphic to each other.

Likewise, all $<$-components of $\conext^\dagger$ are isomorphic to each other.
\end{proposition}

\begin{proof}
Again, we only prove  the claim for $\conext$, as the proof for $\conext^\dagger$ is very similar.  
Let $\tree_0$ and $\tree_1$ be components of $\conext$ and let $\pathp$ be a path in $\tree$. By Lemma~\ref{ch4: paths in F_T}, we can consider a path $\tilde\pathp_0$ in $\tree_0$ and a path $\tilde\pathp_1$ in $\tree_1$ such that $\fftt\big(\tilde\pathp_0\big) = \fftt\big(\tilde\pathp_1\big) = \pathp$. Call $\xi_0$ and $\xi_1$ the sets $\bigcup\{f : f \in \tilde\pathp_0\}$ and $\bigcup\{f : f \in \tilde\pathp_1\}$, respectively. Clearly, $\xi_0$ and $\xi_1$ are functions from $\pathp$ to 2.

Let $f$ be an element of $T_0$ and assume $f \in 2^{\tleq{T}{\ndt}}$. We define the element $\varphi(f)$ of $T_1$ by
\begin{equation}\label{ch4:isomorphic components}
\varphi(f) \in  2^{\tleq{T}{\ndt}} \quad \textrm{and} \quad \varphi(f)(u) := \left\{
\begin{array}{ll}
 f(u) & \textrm{if } u \in  \tleq{T}{\ndt}\setminus \pathp \\
 \xi_1(u) & \textrm{if } u \in   \pathp 
\end{array}\right..
\end{equation}
We prove that $\varphi$ is an isomorphism. It is readily verified that $\varphi$ is order preserving, which implies that it is injective. For every $g \in T_1$, $\varphi^{-1}(g)$ can be defined in the same way as $\varphi(f)$, by exchanging $\tree_0$ and $\tree_1$.
\end{proof}

\paragraph{Observation} Let $\chi$ be any fixed element of $2^T$ and consider the set  $T_\chi = \{ \chi{\upharpoonright}_{\tleq{T}{\ndt}} : \ndt \in T\}$. The correspondence $t \mapsto  \chi{\upharpoonright}_{\tleq{T}{\ndt}}$ is a bijection and $t < t'$ iff $  \chi{\upharpoonright}_{\tleq{T}{\ndt}} \subset  \chi{\upharpoonright}_{\tleq{T}{\ndt'}}$, so that $(T_\chi; \subset)$ is isomorphic to $\tree$. Thus, $\conext$ contains $2^{|T|}$ copies of $\tree$. In this context, we just mention and eventually leave open the question of how many $<$-components $\conext$ has.
Still, as a first step towards answering that question, 
one can observe that the set $\conextset$ can be described as $\{\chi{\upharpoonright}_{\tleq{T}{\ndt}} : \ndt \in T \textrm{ and } \chi \in 2^T\}$. Hence, $f$ and $g$ belong to the same $<$-component of $\conext$ whenever $f = \chi{\upharpoonright}_{\tleq{T}{\ndt}}$, $g = \xi{\upharpoonright}_{\tleq{T}{\ndu}}$, and $\chi$ and $\xi$ coincide on a stem in $\tree$. 
We denote this relationship between $\chi$ and $\xi$ by $\chi \approx \xi$. Then $\approx$ is an equivalence relation and the number of $<$-components of $\conext$ is the number of equivalence classes modulo $\approx$.

\begin{definition}\label{def:condensed extensions}
For every tree $\tree = (T;<)$, a \defstyle{condensed forking} of $\tree$ is a pair $(\tree_c, \fftt_c)$ such that:
\begin{enumerate}
	\item
		$\tree_c = (T_c;<_c)$ is a condensed tree;
	\item
		$\fftt_c$ is an order-preserving function from $T_c$ onto $T$;
	\item
		for every path $\pathp$ in $\tree_c$, the restriction of $\fftt_c$ to $\pathp$ is an isomorphism between $\pathp$ and a path in $\tree$;
	\item
		every path in $\tree$ is the $\fftt_c$-image of a path in $\tree_c$.
\end{enumerate}
If, in addition, 
\begin{enumerate}
	\item[5.]
		there exists a subtree $\tree'_c$ of $\tree_c$ such that the restriction of $\fftt_c$ to $\tree'_c$ is an isomorphism,
\end{enumerate}
then $(\tree_c, \fftt_c)$ will be called a  \defstyle{condensed extension} of $\tree$.
\end{definition} 

For every $<$-component $\tree^*$ of $\conext$, we denote by $\fftt^*$ the restriction of the function $\fftt$ defined in Equation (\ref{ch4:function from F_T to T}) to $T^*$. Then, by Propositions~\ref{ch?:prop:CondensedExt} and \ref{ch4: isomorphic components}, the pair $(\tree^*; \fftt^*)$ is a condensed forking of $\tree$. By the observation above, $(\tree^*; \fftt^*)$ is also a condensed extension of $\tree$.
The same holds likewise for every  $<$-component of $\conext^\dagger$. 

Examples of condensed forkings that are not condensed extensions are the trees $\tree_1^c$ and $\tree_2^c$ in Figure~\ref{fig:ex3:CondensedExtensions} below.

Thus, there can be many non-isomorphic condensed forkings or extensions of a given tree. A natural question arises whether there exists among them a unique one that is smallest by inclusion (up to isomorphism). The following example shows that the answer is negative in both cases.

\begin{example} 
\label{ex3:CondensedExtensions} 
	The trees $\tree^c_1$ and $\tree^c_2$ at the bottom of 
Figure~\ref{fig:ex3:CondensedExtensions} are condensed forkings of $\tree$,  which are smaller in size than  $\conext^\dagger$ and are isomorphically embeddable in it. 
It is easy to see that each of these two smaller condensed forkings of $\tree$ is minimal with this property. In this example they are isomorphic, but they can be made non-isomorphic  by extending the leaves $x,y,z$ in the original tree $\tree$ with pairwise non-isomorphic (and condensed) subtrees. Then, the resulting tree will not have a smallest, up to isomorphism,  
condensed forking.

Similarly, the non-existence of minimal condensed extensions of $\tree$ can be shown by considering the trees obtained by removing $y_1$ or $z_1$ from $\conext^\dagger$.

	\begin{figure}[ht]
		\centering
		\centering
		\begin{picture}(280,170)
			\qbezier[40](30,105)(30,125)(30,135)
			\put(30,105){\circle*{3}}
			\put(30,135){\circle*{3}}
			\qbezier[40](30,135)(20,150)(10,165)
			\put(10,165){\circle*{3}}
			\qbezier[40](30,135)(40,150)(50,165)
			\put(50,165){\circle*{3}}
			\qbezier[40](30,135)(30,155)(30,165)
            \put(30,165){\circle*{3}}
            
			\put(34,103){\footnotesize{$\ndt$}}
			\put(34,130){\footnotesize{$\ndu$}}
			\put(6,170){\footnotesize{$x$}}
			\put(47,170){\footnotesize{$z$}}
			\put(28,170){$y$}
			
			\put(0,115){\footnotesize{$\tree$:}}
			
			
				\qbezier[40](220,105)(205,120)(190,135)
			\put(220,105){\circle*{3}}
                   \put(190,135){\circle*{3}}
                   
			\qbezier[40](190,135)(200,150)(210,165)
			\put(210,165){\circle*{3}}
			\qbezier[40](190,135)(180,150)(170,165)
			\put(170,165){\circle*{3}}
			
			\qbezier[40](220,105)(235,120)(250,135)
			\put(250,135){\circle*{3}}
			\qbezier[40](250,135)(240,150)(230,165)
			\put(230,165){\circle*{3}}
			\qbezier[40](250,135)(260,150)(270,165)
			\put(270,165){\circle*{3}} 
			
			\qbezier[40](190,135)(190,150)(190,165)
			\put(190,165){\circle*{3}} 
			\qbezier[40](250,135)(250,150)(250,165)
			\put(250,165){\circle*{3}} 
	
			\put(224,103){\footnotesize{$\ndt$}} 
			\put(254,130){\footnotesize{$\ndu_1$}}
			\put(195,130){\footnotesize{$\ndu_0$}}
			\put(165,170){\footnotesize{$x_0$}}
			\put(225,170){\footnotesize{$x_1$}}
			\put(204,170){\footnotesize{$z_0$}}
			\put(264,170){\footnotesize{$z_1$}}
			\put(184,170){\footnotesize{$y_0$}}
			\put(244,170){\footnotesize{$y_1$}}
			\put(160,115){\footnotesize{$\conext^\dagger$:}}
		\qbezier[40](30,5)(15,20)(0,35)
			\put(30,5){\circle*{3}}
                   \put(0,35){\circle*{3}}
                   
			\qbezier[40](0,35)(10,50)(20,65)
			\put(20,65){\circle*{3}}
			\qbezier[40](0,35)(-10,50)(-20,65)
			\put(-20,65){\circle*{3}}
			
			\qbezier[40](30,5)(45,20)(60,35)
			\put(60,35){\circle*{3}}
			\qbezier[40](60,35)(50,50)(40,65)
			\put(40,65){\circle*{3}}
			\qbezier[40](60,35)(70,50)(80,65)
			\put(80,65){\circle*{3}} 
	
			\put(34,3){\footnotesize{$\ndt$}} 
			\put(64,30){\footnotesize{$\ndu_1$}}
			\put(5,30){\footnotesize{$\ndu_0$}}
			\put(-25,70){\footnotesize{$x_0$}}
			\put(35,70){\footnotesize{$x_1$}}
			\put(13,70){\footnotesize{$y_0$}}
			\put(73,70){\footnotesize{$z_1$}}
			\put(-25,15){\footnotesize{$\tree^c_1$:}}
		\qbezier[40](220,5)(205,20)(190,35)
			\put(220,5){\circle*{3}}
                   \put(190,35){\circle*{3}}
                 
			\qbezier[40](190,35)(200,50)(210,65)
			\put(210,65){\circle*{3}}
			\qbezier[40](190,35)(180,50)(170,65)
			\put(170,65){\circle*{3}}
			
			\qbezier[40](220,5)(235,20)(250,35)
			\put(250,35){\circle*{3}}
			\qbezier[40](250,35)(240,50)(230,65)
			\put(230,65){\circle*{3}}
			\qbezier[40](250,35)(260,50)(270,65)
			\put(270,65){\circle*{3}} 
			\put(224,3){\footnotesize{$\ndt$}} 
			\put(254,30){\footnotesize{$\ndu_1$}}
			\put(195,30){\footnotesize{$\ndu_0$}}
			\put(165,70){\footnotesize{$x_0$}}
			\put(226,70){\footnotesize{$y_1$}}
			\put(203,70){\footnotesize{$y_0$}}
			\put(263,70){\footnotesize{$z_1$}}
			\put(160,15){\footnotesize{$\tree^c_2$:}}
		\end{picture}

		\caption{A tree $\tree$, the condensed extension $\conext^\dagger$, and 2 smaller (and minimal) condensed forkings of $\tree$ \label{fig:ex3:CondensedExtensions}.}
	\end{figure}
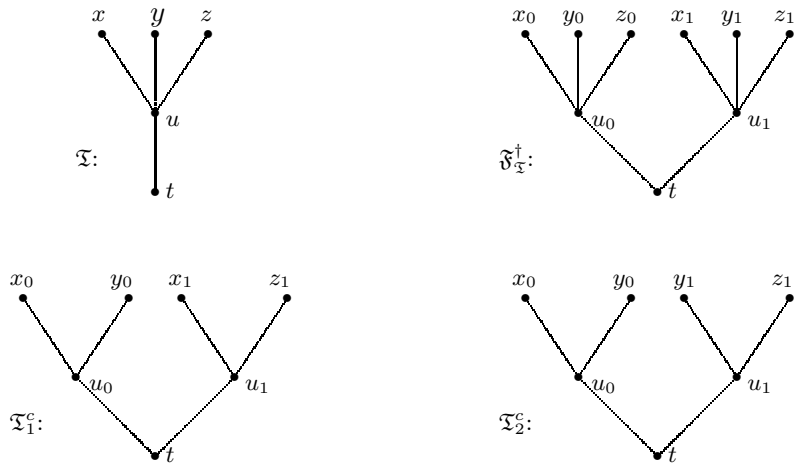
\end{example}

\section{Concluding remarks}
\label{sec:concluding} 

With this work we have initiated an exploration of the general theory of trees. Follow-up research will include: 

\begin{itemize}
\item a study of \emph{complete trees} (in the sense of Dedekind completeness) and constructions of tree completions;

\item a study of general operations on trees, such as sums and products, thus extending classical studies of ordinal arithmetic, due to Cantor, Sierpinski, and others, and, more generally, operations on linear orders (cf. \cite{Rosenstein});

\item a study of classes of trees generated by applying such operations, and their structural and logical theories.
\end{itemize}

Our ultimate goal is a systematic development of a structural theory of trees. One intended target application of this study is to characterise elementary equivalence and other logical equivalences of trees and to obtain new axiomatisations and decidability or undecidability results for logical theories of important classes of trees, in the spirit of those in \cite{GorankoKellerman2021}.

\section*{Acknowledgements}

We thank the referee for the careful reading and helpful comments and suggestions on the paper.

\providecommand{\bysame}{\leavevmode\hbox to3em{\hrulefill}\thinspace}
\providecommand{\MR}{\relax\ifhmode\unskip\space\fi MR }
\providecommand{\MRhref}[2]{%
  \href{http://www.ams.org/mathscinet-getitem?mr=#1}{#2}
}
\providecommand{\href}[2]{#2}

\end{document}